\documentclass[12pt,twoside]{amsart}
\usepackage{amsthm,amsmath,hyperref}
\usepackage{amssymb}
\usepackage{tikz}
\usepackage[all]{xy}
\usepackage{graphicx}
\usepackage{url}
\usepackage{braket}
\usepackage{ifpdf}
\usepackage{mathtools}
\usepackage{color}
\usepackage{rotating}
\usepackage{cleveref}

\newcommand{\mathscr}{\mathcal}
\newcommand{\ba}{\begin{align}}
\newcommand{\ea}{\end{align}}
\newcommand{\be}{\begin{equation}}
\newcommand{\ee}{\end{equation}}

\newcommand{\beq}{\begin{eqnarray}}
\newcommand{\eeq}{\end{eqnarray}}

\newcommand{\beqs}{\begin{eqnarray*}}
\newcommand{\eeqs}{\end{eqnarray*}}

\newcommand\FS{\mathfrak{F}_S}
\newcommand{\GHZ}{|\text{GHZ}\rangle}
\newcommand{\Max}{|\text{Max}\rangle}
\newcommand{\Z}{\mathbb{Z}}
\newcommand{\LL}[1]{{\bf L$_{\bf #1}$}}
\newcommand{\RR}[1]{{\bf R$_{\bf #1}$}}
\newcommand{\Sim}[1]{{\bf S$_{\bf #1}$}}

\tikzstyle WL=[line width=3pt,opacity=1.0]
\hypersetup{
    colorlinks = true,
    linkcolor = blue,
    urlcolor = blue,
    citecolor = blue,
    anchorcolor = black,
}

\begin{document}
\title{A Mathematical Picture Language Program}
\author{Arthur Jaffe}
\author{Zhengwei Liu}

\address{Harvard University, Cambridge, MA 02138, USA}
\address{jaffe@g.harvard.edu}
\address{zhengweiliu@fas.harvard.edu }

\begin{abstract}
We give an overview of our  philosophy of pictures in  mathematics.  
We emphasize a bi-directional process between picture language and mathematical concepts:  abstraction and simulation.  This motivates a program to understand  different subjects, using virtual and real mathematical concepts simulated by pictures.
\end{abstract}

\maketitle

\thispagestyle{empty}

Pictures appear throughout mathematical history, and we recount some  of this story. We explain  insights we gained through using mathematical pictures.
We reevaluate ways that one can use pictures, not only to gain mathematical insights, but also to prove mathematical theorems.
As an example, we describe ways that the quon language, invented to study quantum information, sheds light on several other areas of mathematics.  It results in proofs and new algebraic identities of  interest in several fields. Motivated by this success, we outline a \textit{picture program} for further research.

Our picture language program has the goal to unify ideas from different subjects. We focus here on the three-dimensional quon language~\cite{Quon}.  This language is a topological quantum field theory (TQFT) in 3D space, with lower dimensional defects.   A quon is a 2D defect on the boundary of a 3D manifold. 

We believe quon and other languages will provide a framework for increased mathematical understanding, and 
we expect that looking further into the role of the mathematics of pictures will be productive. Hence in \cref{Sect:Questions}, we  pose a number of problems as the basis for a picture language research program. 

Pictures have been central for visualizing insights and for motivating proofs in many mathematical areas, especially in geometry, topology, algebra, and combinatorics.  They extend from ancient work in the schools of Euclid and Pythagoras to modern ideas in particle physics, category theory, and TQFT.  See the interesting recent account by  \cite{Silver}.
Nevertheless we mention two aspects of mathematical pictures that we feel merit special study.

First is the importance we ascribe to mathematical analysis of pictures.  
We explain how one has begun  to formulate a theory of mathematical analysis on pictures, in addition to the study of their topology and geometry.  For example the analytic aspect of pictures in TQFT is a less-developed area than its topological  and algebraic aspects.  Yet it has great  potential for future advances. 

Second is the notion of proof through pictures.  In focusing on general mathematical properties of pictures, we wish to distinguish this quality from using pictures in a particular concrete mathematical theory.  In other words, we aim to distinguish the notion of the properties of a picture language {\bf L} on the one hand, from its use through a simulation {\bf S} to model a particular mathematical reality {\bf R}.  
We thank one referee for pointing out that the distinction between {\bf L} and {\bf R} parallels the distinction in linguistics between \textit{syntax} and \textit{semantics}.

We propose that it is interesting to prove a  result about the language  {\bf L}, and thereby through simulation ensures  results in {\bf R}. One can use a single picture language {\bf L}  to simulate several different mathematical areas. In fact a theorem in {\bf L} can ensure different theorems in different mathematical subjects \RR1, \RR2, etc., as a consequence of different simulations \Sim1, \Sim2, etc.  This leads to the discussion of the simulation clock in the subsection \textit{What Next?} Different configurations of the hands of the clock reveal the interrelation between picture proofs for seemingly unrelated mathematical results.

We also discuss the important distinction between two types of concepts in {\bf R} that we simulate by a given {\bf S}. These may be {\it real} concepts, or they may be {\it virtual}. This distinction is not absolute, but depends on what language and simulation one considers.  We give some examples, both in mathematics and in physics, in the subsection  {\it Real and Virtual}.
 
We hope these remarks about pictures can enable progress in understanding both mathematics and physics. Perhaps this can even help other subjects, such as  the neurosciences, where ``One picture is worth one thousand symbols.''

\subsection{Euclidean Geometry}
Mathematicians have used pictures since the evolution of Euclidean geometry in ancient Greece.  They proved abstract theorems based on axioms designed from pictorial intuition.
A powerful feature of pictures is that one can easily visualize symmetries. Rotation and reflection symmetries appear in ancient arguments.  

A good example comes from the problem of four points A,B,C,D on a plane, as illustrated in [\ref{FourPoints}].
\vskip -.4cm
\begin{equation}\label{FourPoints}
\raisebox{-1.1cm}{\scalebox{.9}{
\begin{tikzpicture}
\draw (1,0) arc (0:49:1) coordinate (A);
\draw (1,0) arc (0:121:1) coordinate (B);
\draw (1,0) arc (0:169:1) coordinate (C);
\draw (1,0) arc (0:289:1) coordinate (D);
\draw[white] (1.2,0) arc (0:49:1.2) coordinate (A');
\draw[white] (1.2,0) arc (0:121:1.2) coordinate (B');
\draw[white] (1.2,0) arc (0:169:1.2) coordinate (C');
\draw[white] (1.2,0) arc (0:289:1.2) coordinate (D');
\draw (1,0) arc (0:360:1);
\draw (A) to (B);
\draw (B) to (C);
\draw (C) to (D);
\draw (D) to (A);
\draw (A) to (C);
\draw (B) to (D);
\node at (A') {A};
\node at (B') {B};
\node at (C') {C};
\node at (D') {D};
\end{tikzpicture}
}}
\end{equation}
\vskip -.17cm \noindent
The points  lie on a circle if and only if the angles BAC and BDC are equal. This can be established pictorially  or algebraically, and the picture 
proof\footnote{{
The argument depends on the \textit{inscribed angle theorem}: Three points B, D, C on a circle determine the angle BDC $=\theta$, and the  angle BOC $=\psi$, where O is the center of the circle. Then always $2\theta=\psi$.
The proof: In the special case that BD passes  through O,
symmetry shows that the  triangle COD is  isosceles. As the sum of angles in a triangle is $\pi$,  both $2\theta$ and $\psi$ complement the same angle, so  $2\theta=\psi$. The general case then follows  by drawing the diameter through BO, and considering the sum or difference of two special cases.}}
is elementary.

\textheight = 8.6in 

 \subsection{Abstraction and Simulation} 
 We propose two basic components to understanding that we call ~{\bf L} and~{\bf R}.  Here~{\bf L}  denotes \textit{abstract concepts} or language, while  {\bf R} stands for the  \textit{concrete subjects} or reality, which we desire to understand.  We could also think of them as \textit{left} and \textit{right}.
Simulation {\bf S} represents a map from {\bf L} to {\bf R}.

Our universe provides a great reservoir for ideas about the real world. We can consider this as~{\bf R}. We understand these ideas through abstraction, including theories of  mathematics, physics, chemistry, and biology.  To deal with these real concepts, one often requires 
virtual concepts that have no meaning in the real world.  
These virtual concepts may not have an immediate real meaning, yet they may provide  key insight to understanding the real structures. 
\goodbreak

Abstraction is a method  to study complicated ideas and goes from~{\bf R}  to~{\bf L}. 
For example, chemistry as~{\bf L} may provide a logical language to abstract certain laws in biology as~{\bf R}, with a simulation {\bf S: L $\rightarrow$ R}. But one can continue this chain where we regard chemistry as a new~{\bf R} and physics as its abstraction as a new~{\bf L} with a new simulation~{\bf S}.  Abstraction can be repeated yet again with mathematics as~{\bf L} and physics as~{\bf R}.  
At each step one learns the axioms in~{\bf L} from the real concepts in~{\bf R}.  In order to do computation, we often need to introduce virtual concepts in~{\bf L} to understand the concepts in~{\bf R}.

We are especially interested in the case that {\bf L} is a picture language.
In picture language, one can represent the concepts in {\bf L}  by pictures which play the role of logical words, with axioms as their grammar.
The axioms should be compatible with pictorial intuition.
One can develop a picture language as an independent theory like Euclidean geometry. 

A good mathematical simulation for a particular mathematical subject {\bf R}, should satisfy two conditions.  Mathematical concepts in {\bf R} should be simulated by simple pictures in {\bf L},  and mathematical identities in {\bf R} should arise from  performing elementary operations 
on the pictures in {\bf L}.
Abstraction and simulation could be considered as inverse processes of each other.

\subsection{Real and Virtual} \label{Sect:Real-Virtual}
The idea of whether a concept is real or virtual plays an important role.  This is not an absolute concept, but depends on the simulation.  For example, over the real-number field, the square root $\sqrt{x}$  of a positive number $x$ is real, while the square root of a negative number $x$ is virtual.  By enlarging the field  to complex  numbers, one  understands this virtual concept in a fruitful way.  
One gets used to the power of complex numbers and in a new simulation,  one can regard them as ``real.''

In a picture language, it may be that a picture can represent a real concept in one theory and a virtual one in another. As an example consider Feynman diagrams, the pictures that one uses  to describe particle interactions. A contribution to the scattering of two physical electrons is described by the first electron  emitting a photon (quantum), that is absorbed by the second electron.  However, conservation of energy and momentum preclude the exchanged photon from being real: it must be virtual.  Fig. 1 shows the original diagram in~\cite{Feynman}.  

Analytic continuation (Wick rotation) to imaginary time puts the real and virtual particles on equal footing. In that situation the above example only involves virtual particles. In this way, one will obtain many virtual concepts, which could become real based on a proper simulation.

\begin{figure}
\begin{center}
\includegraphics[scale=.4]{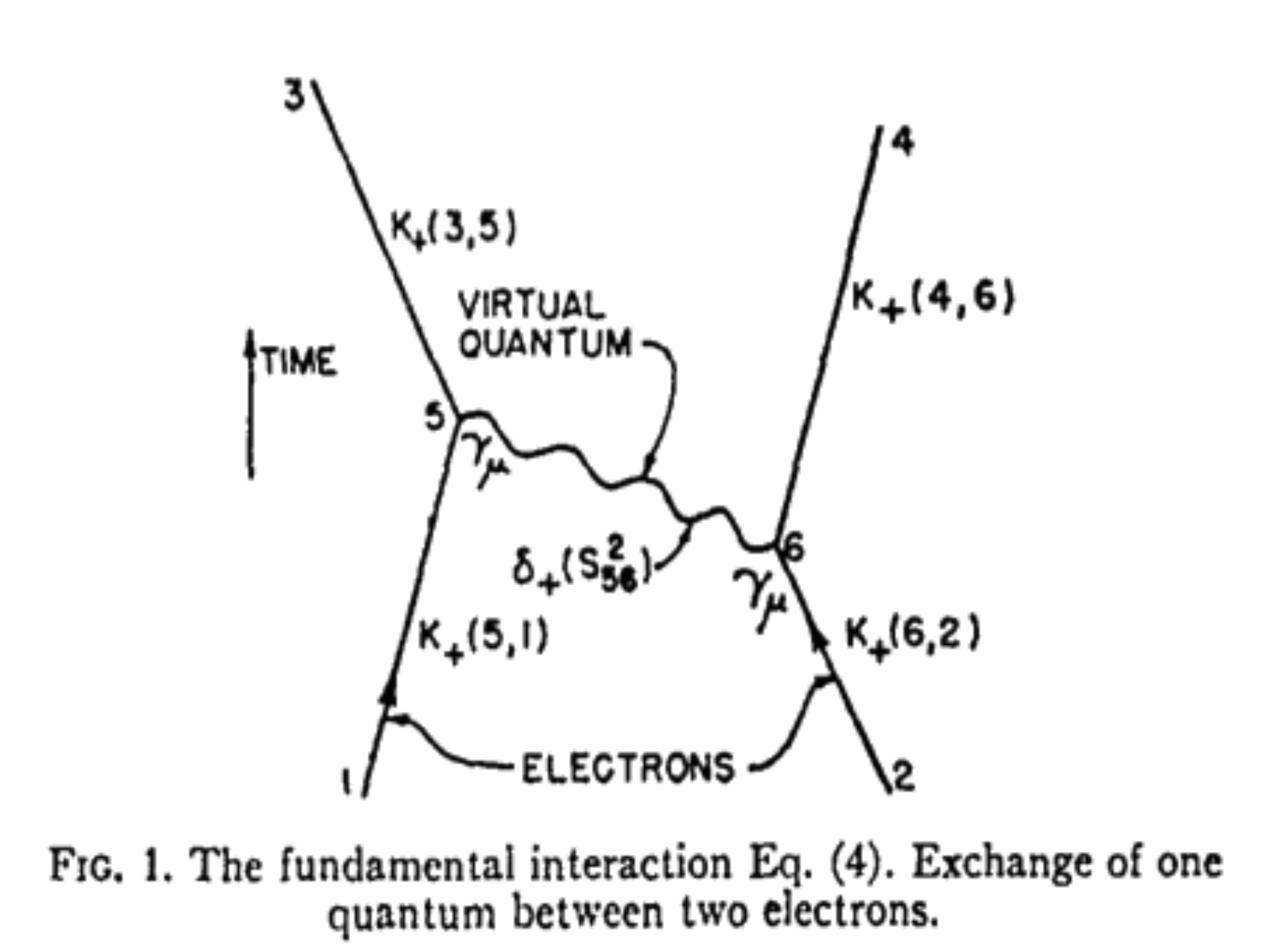}
\end{center}
\vskip -.6cm
\caption{The First Feynman Diagram, reprinted with APS permission  from \cite{Feynman}.}
\end{figure}
\subsection{Modern Picture Mathematics}
Pictures played a central role both in the invention of topology, as well as in its modern understanding.
A synergetic algebraic formalization was developed in parallel in category theory, which originated in the work of Eilenberg and MacLane in the 1940's~\cite{Eilenberg-MacLane}.

TQFT appeared as a way to understand various topics using cobordisms.  Remarkable examples of this point of view can be found in the work of Jones \cite{Jones}, Witten \cite{Witten}, Atiyah~\cite{Atiyah}, Reshetikhin-Turaev~\cite{RT}, Turaev-Viro~\cite{Turaev-Viro}, and Ocneanu~\cite{Ocneanu}.  Various generalizations inspired by TQFT have been studied. And many results in this direction have emerged over the past thirty years. 

Compared with topology, it is less obvious that pictures can also shed light on the study of algebra---in particular on representation theory.  It turns out that this is the case.  The insights from pictures are   useful technically (such as with  Young diagrams,  braids, or  quivers).  These insights are also useful conceptually (e.g. to show how the topological properties of multi-particle systems are captured by properties of the centralizer algebra of representations, in the sense of Schur-Weyl duality \cite{Schur,Weyl}).

Last, but not least, pictures also have deep connections with analysis on infinite-dimensional Hilbert space. This led to the discovery of the Jones polynomial~\cite{JonesPoly}.  Conformal field theory (CFT) is a topic {\bf R} from which we can learn laws for a picture language {\bf L}.  

\section{Beyond Topology}
When Atiyah defined TQFT mathematically, he wrote 
``... it may well be that such topological understanding is a necessary pre-requisite to building the analytical apparatus of the quantum theory''~\cite{Atiyah}.
Today we want to move forward to reach a full quantum field theory. Our long-term goal is to construct quantum field theory using pictures. However, that may be too hard to achieve  as the initial step.

We emphasize two properties of QFT that shed light on pictures. The first theme is \textit{symmetry}.
In order to understand continuum symmetry, it is useful to understand discrete symmetry that can approximate the continuum.

The second theme \textit{positivity} is especially interesting,  as positivity provides the basis for  analysis.  The striking fact is that the analysis of pictures is a mathematical theory with great potential, but still in its infancy. On one hand, we expect to abstract  concepts from analysis to pictures. This will enrich the theory of picture language in one additional dimension. On the other hand, we want to simulate analysis using picture language and provide new pictorial tools.

\subsection{Symmetry}
An elementary characteristic that transcends topology is \textit{shape}.
This goes in the direction of encapsulating geometry of the picture.  The geometry can indicate the presence of additional \textit{symmetry} represented by the pictures.  With the shape of a picture as an additional tool, one can ask where it leads.

In lattice models of statistical physics, people often study square lattices and honeycomb lattices, which capture additional symmetry in two and three directions respectively. 
The pictorial duality of lattices could provide interesting dualities of mathematical theories, such as the analog of Kramers-Wannier duality~\cite{Kramers-Wannier}, illustrated in [\ref{LatticeDuality}],
\begin{equation}\label{LatticeDuality}
\raisebox{-.8cm}{
\scalebox{.3}{
\begin{tikzpicture}
\foreach \x in {0,1,2,3}
{
\draw (\x,-1)--+(0,5);
\foreach \y in {0,1,2,3}
{
\draw (-1,\y)--+(5,0);
\node at (\x,\y) {$\bullet$};
}}
\fill[white] (-1,-1) rectangle (-.8,4);
\fill[white] (4,-1) rectangle (3.8,4);
\fill[white] (-1,-1) rectangle (4,-.8);
\fill[white] (-1,4) rectangle (4,3.8);
\begin{scope}[shift={(6.3,1.5)},xscale=.8,yscale=.8]
\node at (-1,0) {$\longrightarrow$};
\end{scope}
\begin{scope}[shift={(8,0)}]
\foreach \x in {-1,0,1,2,3,4}
{
\foreach \y in {-1,0,1,2,3,4}
{
\fill[gray!50]  (\x-.5,\y)--(\x,\y-.5)--(\x+.5,\y)--(\x,\y+.5);
}}
\foreach \x in {0,1,2,3}
{
\draw (-1,-.5+\x) --++ (4.5-\x,4.5-\x);
\draw (\x-.5,-1) --++ (4.5-\x,4.5-\x);
\draw (-1,3.5-\x) --++ (4.5-\x,-4.5+\x);
\draw (\x-.5,4) --++ (4.5-\x,-4.5+\x);
} 
\fill[white] (-1.5,-1) rectangle (-.8,4.5);
\fill[white] (4.5,-1) rectangle (3.8,4.5);
\fill[white] (-1.5,-1.5) rectangle (4.5,-.8);
\fill[white] (-1.5,4.5) rectangle (4.5,3.8);
\end{scope}
\end{tikzpicture}
}}\;.
\end{equation}
If a vector space $V$ is simulated by square-like pictures on a two dimensional plane, then gluing two pictures vertically or gluing them horizontally, defines two multiplications of $V$. 
The $90^{\circ}$  rotation,  called the string Fourier transform (SFT) and that we denote by  $\FS$, intertwines the two multiplications as illustrated below.  See \cite{PPA,Quon, SFT-LW} for details and further references.
\begin{center}
\hskip -.3in Fourier Transform \quad \quad Multiplication \quad \quad Convolution 
\end{center}
\begin{equation*}
\begin{tikzpicture}
\hskip -.2in
\begin{scope}[shift={(-2,1.2)},scale=.25]
\fill[blue!20] (0,0) rectangle (1,1);
\draw (0,0) rectangle (1,1);
\draw (0,0)--(-1,2);
\draw (0,1)--(2,2);
\draw (1,1)--(2,-1);
\draw (1,0)--(-1,-1);
\node at (3,.2) {$:$};
\end{scope}

\hskip .4in
\begin{scope}[scale=.25]
\foreach \x in {0,1}{
\foreach \y in {0,1}{
\foreach \u in {0}{
\foreach \v in {1,2,3}{
\coordinate (A\u\v\x\y) at (\x+1.5*\u,\y+3*\v);
}}}}

\foreach \u in {0}{
\foreach \v in {1,2}{
\fill[blue!20] (A\u\v00) rectangle (A\u\v11);
\draw (A\u\v00) rectangle (A\u\v11);
}}

\draw (0,3)--++(-.5,-.5);
\draw (1,3)--++(.5,-.5);
\draw (0,7)--++(-.5,.5);
\draw (1,7)--++(.5,.5);

\draw (A0101) to [bend left=30] (A0200);

\draw (A0111) to [bend left=-30] (A0210);

\node at (5,5) {$\rightarrow$};

\begin{scope}[shift={(16,4.5)},rotate=90]
\foreach \x in {0,1}{
\foreach \y in {0,1}{
\foreach \u in {0}{
\foreach \v in {1,2,3}{
\coordinate (A\u\v\x\y) at (\x+1.5*\u,\y+3*\v);
}}}}

\foreach \u in {0}{
\foreach \v in {1,2}{
\fill[blue!20] (A\u\v00) rectangle (A\u\v11);
\draw (A\u\v00) rectangle (A\u\v11);
}}

\draw (0,3)--++(-.5,-.5);
\draw (1,3)--++(.5,-.5);
\draw (0,7)--++(-.5,.5);
\draw (1,7)--++(.5,.5);

\draw (A0101) to [bend left=30] (A0200);

\draw (A0111) to [bend left=-30] (A0210);
\end{scope}
\end{scope}
\end{tikzpicture}
\end{equation*}
These 2D pictorial operations coincide with Fourier transform, multiplication, and  convolution in Fourier analysis. This is a cornerstone for understanding pictorial Fourier duality.

\subsection{Analysis}
How do we connect pictures with analysis?
Usually pictures without boundary are scalars. How can we formulate  a measurement pictorially?
A main lesson from quantum field theory is the importance of reflection positivity. An elementary pictorial interpretation of reflection positivity  is that gluing a picture to its mirror image is positive.

In constructive quantum field theory, extensive analysis is performed through estimating Feynman diagram pictures in terms of subdiagrams; this may use a pictorial Schwarz inequality, or some more sophisticated operator norm. These estimates are central to the proof of most results in the subject.  
In the framework of our present discussion, such analysis takes place on pictures in {\bf R}. 
We call this using \textit{pictures in analysis}. 

But we are really interested in whether, and to what extent, one can do \textit{analysis on pictures}.  This means that we need to do computations in {\bf L}, without reference to {\bf R}.  
The discussion of rotation, Fourier analysis, multiplication, and convolution indicate that some progress can be made. 

Do we obtain interesting analysis based on this minimal requirement? In fact the surprise is that we already have interesting results inspired by Fourier theory.
As mentioned, the 2D pictorial operation on square-like pictures is compatible with Fourier analysis.

Compactifying the plane to a sphere defines a measurement based on reflection positivity. Pictorial consistency on the sphere implies that the Fourier transform is a unitary. Many other results in Fourier analysis carry over to pictures.

\section{A Pictorial Journey}
Freeman Dyson described two types of mathematicians: birds and frogs~\cite{Dyson}. They do mathematics in different ways. Birds soar between different fields with unifying ideas, like Yuri Manin in his book \textit{Mathematics as Metaphor}~\cite{Manin}.  Frogs sort out details to achieve great depth of understanding.
In fact it is good to attempt to encompass both metaphors!  One does  find both these ingredients in our story of picture language.  While we began in quantum information, we ended up traveling through much of the colorful landscape of mathematics.

In our context of language, the bird flies back and forth to discover new {\bf R}'s, {\bf L}'s, and {\bf S}'s.  The frog uses an {\bf S} to understand some important problem.
The shape of pictures provides key hints and insights for finding these connections. A popular article by Peter Reuell described picture language as lego-like mathematics \cite{Reuell}.

\subsection{Harvard}
The authors began our collaboration two years ago in late July 2015  with much discussion at Harvard. Our first goal was to understand reflection positivity for parafermions  \cite{RPPF,CRP} in a pictorial way~\cite{PPA}.  If we call this mathematical problem  \RR1, then it led us to define the picture language  \LL1 which we call planar para algebra.  This is a generalization of planar algebra, but with strings replaced by charged strings, and topological isotopy replaced by para isotopy. The map from \LL1 to \RR1 we call the simulation \Sim1.

In the simulation  \Sim1, we found an elementary explanation that a $90^{\circ}$ rotation of pictures represented Fourier transformation, and takes multiplication to convolution.
We used this fact to give a geometric, picture proof of reflection positivity~\cite{PPA}. 
\vskip -.3cm
\begin{equation*}
\begin{tikzpicture}
\begin{scope}[shift={(0,-.5)},xscale=.5]
\node at (1,2) {Reflection Positivity};
\foreach \x in {-1,0,1,2}{
\foreach \y in {-1,0,1,2}{
\foreach \z in {-1,0,1,2}{
\coordinate (A\x\y\z) at (\x+\y*.4,\z+\y*.5);
\coordinate (B\x\y\z) at (1+\x+\y*.4,\z+\y*.5);
}}}

\fill[blue!50] (A000)--(A100)--(A110)--(A111)--(A011)--(A001)--(A000);

\fill[red!50] (B000)--(B100)--(B110)--(B111)--(B011)--(B001)--(B000);

\foreach \y in {0,1}{
\draw (A0\y0)--(A1\y0)--(A1\y1)--(A0\y1)--(A0\y0);
}
\foreach \x in {0,1}{
\foreach \z in {0,1}{
\draw (A\x0\z)--(A\x1\z);
}}

\foreach \y in {0,1}{
\draw (B0\y0)--(B1\y0)--(B1\y1)--(B0\y1)--(B0\y0);
}

\foreach \x in {0,1}{
\foreach \z in {0,1}{
\draw (B\x0\z)--(B\x1\z);
}}
\end{scope}

\node at (2,0.3) {$\xrightarrow{\phantom{SFT}}$};

\begin{scope}[shift={(3,0.1)},yscale=.5]
\node at (.8,2.8) {$C^*$ Positivity};

\foreach \x in {-1,0,1,2}{
\foreach \y in {-1,0,1,2}{
\foreach \z in {-1,0,1,2}{
\coordinate (A\x\y\z) at (\x+\y*.4,\z+\y*.5);
\coordinate (B\x\y\z) at (\x+\y*.4,-1+\z+\y*.5);
}}}

\fill[red!50] (B000)--(B100)--(B110)--(B111)--(B011)--(B001)--(B000);
\fill[blue!50] (A000)--(A100)--(A110)--(A111)--(A011)--(A001)--(A000);

\foreach \y in {0,1}{
\draw (A0\y0)--(A1\y0)--(A1\y1)--(A0\y1)--(A0\y0);
}
\foreach \x in {0,1}{
\foreach \z in {0,1}{
\draw (A\x0\z)--(A\x1\z);
}}

\foreach \y in {0,1}{
\draw (B0\y0)--(B1\y0)--(B1\y1)--(B0\y1)--(B0\y0);
}

\foreach \x in {0,1}{
\foreach \z in {0,1}{
\draw (B\x0\z)--(B\x1\z);
}}

\end{scope}

\end{tikzpicture}
\end{equation*}
We also found elementary picture-representations for $d\times d$ unitary Pauli matrices $X,Y,Z$ with eigenvalues $q^{j}$, where   $q=e^{2\pi i/d}$, and $j=0,1,\ldots,d-1\in \mathbb{Z}_{d}$.  

Alex Wozniakowski pointed out that our work seemed related to quantum information.  This led to fruitful collaboration among the three of us, in which we used the language \LL1  to simulate quantum communication in the mathematical framework that we name  \RR2. Using this simulation \Sim2 
 one picture deformed by isotopy into different shapes simulates different concepts in quantum information.  In this way we reproduced the standard teleportation protocol of Bennett et al \cite{Bennett} by a topological design in \LL1~\cite{Holographic}.  With these concepts in place, we could  also design  other protocols, including new multi-partite, teleportation protocols~\cite{ConstructiveSimulation}.

The important point for our understanding of  \RR2 was to follow pictorial intuition. This led to finding a natural candidate for a resource state that we called $\Max$.  Our picture Max in \LL1 for  the entangled state $\Max$ is simple and natural, as well as suggesting entanglement in a pictorial fashion. Our 
2-string picture for $\Max_{n}$ (with $n$ qudits) is 
\begin{equation*}
\raisebox{-.7cm}{Max$_{n}={\frac{1}{d^{n/4}}}$}
\scalebox{.7}{
$\underbrace{
\raisebox{-1.5cm}{
\tikz{
\foreach \x in {0,2}{
\draw (\x,0) arc (180:0:.5);
}
\node at (1.5,0) {$\cdots$};
\draw (-1,0) arc (180:0:2.5 and 1.5);
}}}_{2n}$
}
\end{equation*}
The algebraic formula for the simulation  \Sim2{Max}$_{n}$ in \RR2 is 
\be
\Max_{n}={\frac{1}{d^{(n-1)/2}}}\sum_{|\vec{k}|=0} \ket{\vec{k}}_{n},
\ee
where we call $|\vec{k}|=k_{1}+\cdots+k _{n}$ the total charge in $\Z_{d}$.
It involves $d^{n-1}$ terms for a resource state with  $n$ qudits, so algebraically the simulation of Max$_{n}$   is complicated, see~\cite{Holographic}. 

\subsection{ETH Zurich}
We spent four months visiting the Research Institute for Mathematics of the ETH, and finally had time to begin writing up these results~\cite{PPA} and later results~\cite{Holographic,ConstructiveSimulation}.
We then learned that  Greenberger, Horne, and Zeilinger had long before introduced another  multi-qubit resource state in quantum information. This $\GHZ$ state appears to be algebraically simpler, as it is a sum of only $d$ terms,  independent of the number $n$ of qubits,
\be
\GHZ_{n}={\frac{1}{d^{1/2}}}\sum_{k\in Z_{d}} \ket{k,k,\ldots, k}_{n}.
\ee
Both $\Max$ and $\GHZ$ generalize the Bell states, and they have some very similar properties.  This led us eventually to observe that $\GHZ$ and $\Max$ are related by a change of basis---in fact by Fourier transformation on $\mathbb{Z}_{d}$,
	\be
	\Max_{n} = F^{\otimes n} \GHZ_{n}\;.
	\ee

In order to understand this further, we were intrigued by our generalization of ``Kitaev's map'' from Majoranas to the Pauli spin matrices $X,Y,Z$.  We  had discovered a natural generalization to represent  a single qudit as a neutral parafermion/anti-parafermion pair~\cite{PPA}.  Neutrality provides an elegant way to reduce the $d^{2}$ dimensional space of states for two virtual parafermions to the correct $d$  dimensional space for one real qudit.  It involved introducing a new ``four-string'' planar language \LL2 to describe a single qudit, and a corresponding simulation \Sim3\LL2$=$\RR3.  The model \RR3 contains $d$ real and $d^{2}-d$ virtual one-qudit states. The representations of the qudit Pauli $X,Y,Z$ matrices are neutral, so they act on the neutral (real) subspace of $d$ dimensions.  Another hint that neutrality holds the key, is that the SFT equals the discrete Fourier transform in the neutral subspace of  \LL2.

But the language \LL2 and simulation \Sim3 pose a problem for describing more than one qudit.  Braiding charged strings from different qudits destroys neutrality of the individual excitations.  This would allow transitions from the real $n$-qudit space of dimension $d^{n}$, into the virtual $n$-qudit space of dimension $d^{2n}$. So the obstruction to describing multi-qudit states boiled down to the question: how can  one ensure that real multi-qudit states evolve into real multi-qudit states? When we met Daniel Loss in Basel, we found that he was also considering this question.  We did not find the answer immediately. 

\subsection{Bonn}
After Zurich, we had the opportunity to spend six weeks visiting two institutes in Bonn. 
 We discovered the answer to the puzzle described above during June 2016, perhaps receiving inspiration from working in Fritz Hirzebruch's former office at the Max Planck Institute for Mathematics.

During that time, we encountered two more hints about \LL2.  The first is that the Frobenius algebra for the $m$-interval Jones-Wassermann subfactor  in CFT~\cite{Xu} is 
\be
\gamma=\bigoplus_{\vec{X} }\dim(\vec{X})\,\vec{X}.
\ee
Here $\vec X = X_{1}\otimes \cdots \otimes X_{m}$, a tensor of simple objects in a modular tensor categories (MTC) \cite{RT}, and $\dim(\vec{X})$ is the multiplicity of 1 in $\vec{X}$.
This formula coincides with $\Max$ for the group $Z_{d}$, but they have completely different meanings. Secondly, we learned the relations for  bi-Frobenius algebras in the manuscript for~\cite{CoeckeKissinger}, which one of the authors had shared with Alex Wozniakowski. 

The construction of Jones-Wassermann subfactors for an MTC requires an extension of pictures from 2D to 3D space. Here the coincidence between $\Max$ and $\gamma$ is explained through a new $m-n$ duality~\cite{Liu-Xu}.
Going from 2D to 3D also gives a natural explanation of the  bi-Frobenius algebras relation. 

The last piece in solving the puzzle came by adding three manifolds to these 3D pictures; this was inspired by TQFT. Finally we unified all these ideas by formulating  the 3D quon language \LL3, and a simulation \Sim 4 to quantum information. Moreover, we extended our approach from $Z_{d}$ symmetry to MTCs based on work about the Jones-Wassermann subfactor.

We designed \LL3 using ideas from different {\bf R}s: quantum information, subfactor theory, TQFT, and CFT.  Therefore we expected to simulate those  {\bf R}s using \LL3 as well as others.

In the quon language, the bi-Frobenius algebra relation has a topological interpretation.  Both $\Max$ and $\GHZ$ are represented  by single pictures Max and GHZ, where one picture is a $90^{\circ}$ rotation of the other.  Algebraically one resource state is the Fourier transform of the other.
Moreover, the complicated resource state $\Max$ can be computed by using its relation to the elementary resource state $\GHZ$. For 3-qudits the quon state pictures are just rotations of one another, as illustrated in [\ref{GHZ to MAX}]:
\begin{equation}\label{GHZ to MAX}
\scalebox{1}{
\raisebox{-.9cm}{
\begin{tikzpicture}
\node at (-1.7,1.5) {GHZ$_{3}=$};
\node at (2,1.5) {Max$_{3}=$};
\begin{scope}[scale=.22]
\draw[blue] (.5,6.5) circle (4);
\foreach \x in {0,1}{
\foreach \y in {0,1}{
\foreach \u in {0}{
\foreach \v in {1,2,3}{
\coordinate (A\u\v\x\y) at (\x+1.5*\u,\y+3*\v);
}}}}

\foreach \u in {0}{
\foreach \v in {1,2,3}{
\fill[blue!20] (A\u\v00) rectangle (A\u\v11);
\draw (A\u\v00) rectangle (A\u\v11);
\draw[->] (A\u\v00)--++(0,.5);
}}

\draw (A0101) to [bend left=30] (A0200);
\draw (A0201) to [bend left=30] (A0300);
\draw (A0301) to [bend left=-30] (A0100);

\draw (A0111) to [bend left=-30] (A0210);
\draw (A0211) to [bend left=-30] (A0310);
\draw (A0311) to [bend left=30] (A0110);

\begin{scope}[shift={(24,6)},rotate=90]
\draw[blue] (.5,6.5) circle (4);
\foreach \x in {0,1}{
\foreach \y in {0,1}{
\foreach \u in {0}{
\foreach \v in {1,2,3}{
\coordinate (A\u\v\x\y) at (\x+1.5*\u,\y+3*\v);
}}}}

\foreach \u in {0}{
\foreach \v in {1,2,3}{
\fill[blue!20] (A\u\v00) rectangle (A\u\v11);
\draw (A\u\v00) rectangle (A\u\v11);
\draw[->] (A\u\v01)--++(.5,0);
}}

\draw (A0101) to [bend left=30] (A0200);
\draw (A0201) to [bend left=30] (A0300);
\draw (A0301) to [bend left=-30] (A0100);

\draw (A0111) to [bend left=-30] (A0210);
\draw (A0211) to [bend left=-30] (A0310);
\draw (A0311) to [bend left=30] (A0110);
\end{scope}
\end{scope}
\end{tikzpicture}}
}
\end{equation}
Furthermore, the quon language expressed elegantly in pictures the duality between orthonormal bases for the Pauli matrices $X$ and $Z$, so one could begin to imagine pictures describing quantum coordinates.

\subsection{Back at Harvard}
\begin{figure}[h]
\begin{center}
\includegraphics [scale=.25] {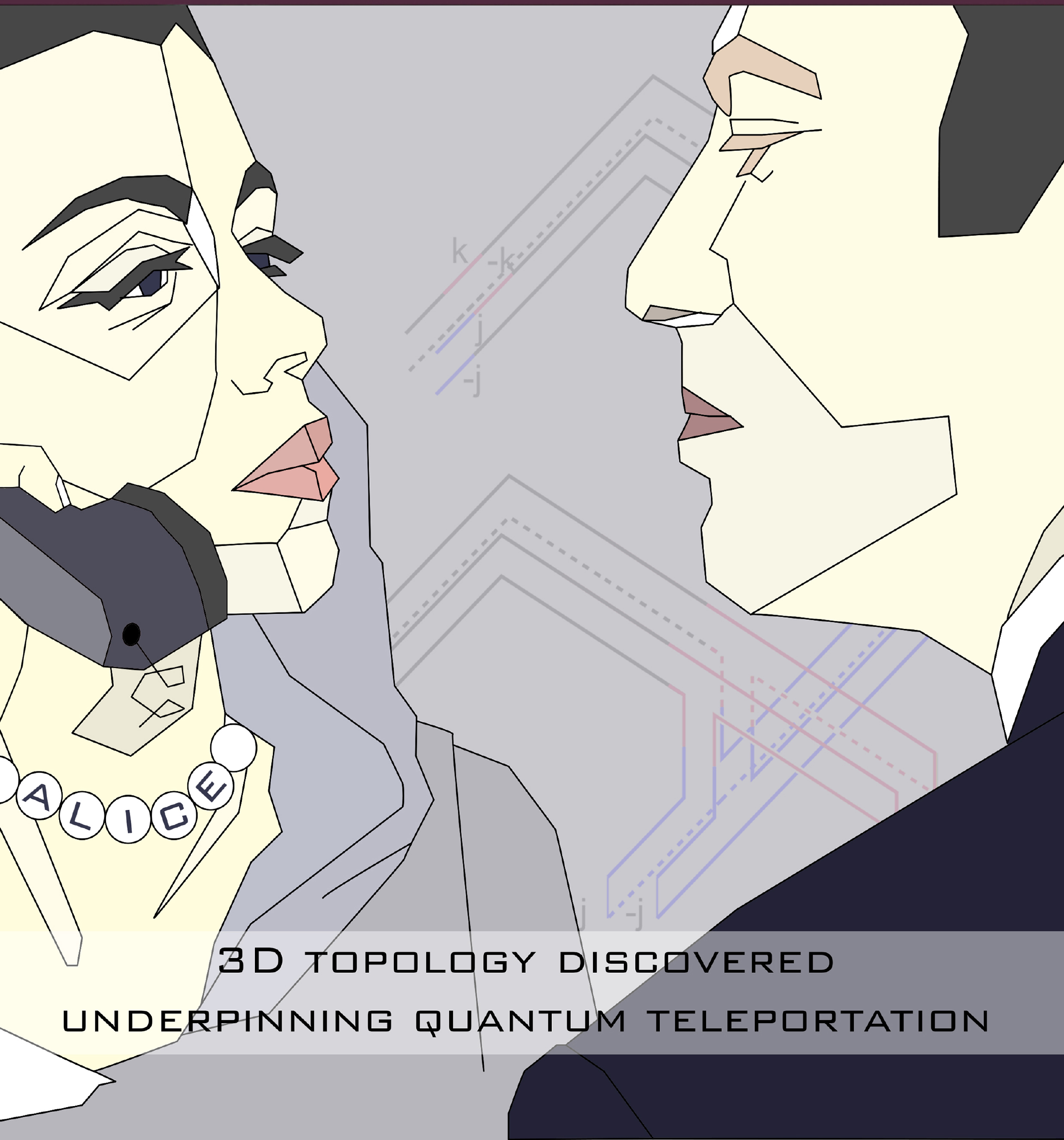}
\end{center}
\caption{Three-Dimensional Representation of Quantum Teleportation}
\end{figure}
Our first discovery by \Sim 4 was  the topological nature of the Feynman, or controlled NOT (CNOT), gate.    The above drawing by Lusa Zheglova captures our representation of the quantum teleportation protocol.  It illustrates classical communication in the foreground. In the background it shows the Quon representation of a Bell state, a CNOT gate, the Fourier transform,  a measurement, and Pauli matrices (corresponding to Kitaev’s map) used in the recovery map. Together they give the widely used quantum teleportation protocol of Bennett et al., and illustrate its elegant 3D topological interpretation.   The details can be found in~\cite{Quon}.

This new quon language \LL3 provided an opportunity to collect ideas and to think about their meaning.  It became apparent that the quon language \LL3 was important on its own, as it could have applications in other areas of mathematics and physics  besides quantum information.  Furthermore the algebraic identity between $\GHZ$ and $\Max$ given by the rotation of the diagrams could provide insight  in other subjects through use of other simulations.  

So what is the meaning of this interesting algebraic identity?  
It actually is the Verlinde formula.
Mathematically, one can generalize the pictorial construction to define GHZ and Max on higher genus surfaces. 
Then the pictorial Fourier duality between GHZ and Max leads to the generalized Verlinde formula~\cite{Verlinde}  for any MTC on the genus $g$ surface,
\begin{align*}
\text{Max}_{n,g}&=\FS^{\otimes n} \text{GHZ}_{n,g}\;,
\\
\Longrightarrow \quad
\dim(\vec{k},g)&=\sum_{k}\left(\prod_{i=1}^n S_{k_i,k} \right)S_{k,0}^{2-n-2g} \;.
\end{align*}
Here on the first line $n$ is the number of quons, 
$\FS$ is the string Fourier transform, and on the second line $\vec{k}=(k_{1},k_{2},\ldots,k_{n})$ are punctures (or marked points) on the genus-$g$ surface, $S$ is the modular transformation, dim is the dimension of the associated (moduli) space.

This pictorial Fourier duality also coincides with the duality of graphs on the sphere.
The dual graph of the tetrahedron is also a tetrahedron. Applying the quon language to this graphic duality, we obtain a general algebraic identity \eqref{6jRelation} for $6j$-symbol self duality.   With $\overline{X}$ denoting the dual object to $X$ in an MTC, 
\be\label{6jRelation}
\left|{{{X_{6}~X_{5}~X_{4}}\choose{\overline{X_{3}}~\overline{X_{2}}~\overline{X_{1}}}}}\right|^{2}
= \sum_{\vec {Y}} \left(\prod_{k=1}^{6}S_{X_{k}}^{Y_{k}} \right)
\left|{{{Y_{1}Y_{2}Y_{3}}\choose{Y_{4}Y_{5}Y_{6}}}}\right|^{2}\;.
\ee
In the special case of quantum $SU(2)$, this was discovered by Barrett~\cite{Barrett}, based on an interesting  identity of  Roberts~\cite{Roberts}.  
The general formulation and proof  of \eqref{6jRelation} is in \S{6} of~\cite{Liu}.

For each graph on a surface, the graphic duality gives a new algebraic identity for MTCs in quon language. 
Most of these identities have virtual meanings. 
Each graph can also be considered as a linear functional generalizing integration.
The pictures generalize the symbol $\int$ and capture additional pictorial relations, such as graph duality mentioned above. 
It would be interesting to understand these new identities and integrations in some new {\bf R}.

\smallskip
\noindent {\bf What Next?}
The progression 
	\begin{center}
	\RR1 $\rightarrow$ \LL1 $\rightarrow$
	\RR2 $\rightarrow$ \LL2 $\rightarrow$
	\RR3 $\rightarrow$ \LL3 $\rightarrow$
	\RR4,
	\end{center}
leads one to believe that much more is in store for the future.  Each progression was inspired by insight from the previous step. We expect that this sequence will continue. We provide some possible {\bf R}s in a simulation clock! 
\vskip -1cm
\begin{center}
Simulation Clock
\end{center}
\begin{equation*}
\scalebox{.5}{
\begin{tikzpicture}
\begin{scope}[scale=1]
\coordinate (O) at (0,0);
\foreach \x in {1,2,3,4,5,6,7,8,9,10,11,12}
{
\draw[white] (0,5.5) arc (90:90-360/12*\x:5.5) coordinate (A\x);
\draw[white] (0,5) arc (90:90-360/12*\x:5) coordinate (AA\x);
\draw[white] (0,2) arc (90:90-360/12*\x:2) coordinate (B\x);
\draw[white] (0,3.3) arc (90:90-360/12*\x:3.3) coordinate (C\x);
\foreach \y in {1,2,3,4}{
\draw[white] (0,5.3) arc (90:90-360/12*\x -6*\y:5.3) coordinate (D\x\y);
\draw[white] (0,5.7) arc (90:90-360/12*\x -6*\y:5.7) coordinate (E\x\y);
}}
\foreach \x in {1,2,3,4,5,6,7,8,9,10,11,12}
{
\node at (AA\x)  {\scalebox{3}{\textcolor{blue!50}{\x}}};
\foreach \y in {1,2,3,4}{
\draw (D\x\y)--(E\x\y);
}}
\fill[blue] (O) circle (.2);
\draw[white] (0,3) arc (90:-6+90-360/12*10:3) coordinate (H);
\draw[white] (0,1.5) arc (90:-6+90-360/12*10:1.5) coordinate (H');
\draw[white] (0,4.5) arc (90:90-360/12*2:4.5) coordinate (M);
\draw[white] (0,2.5) arc (90:90-360/12*2:2.5) coordinate (M');
\draw[thick,->,blue] (O)--(H);
\draw[->,blue] (O)--(M); 
\node at (H') {\rotatebox{-215}{\text{Quons}}};
\node at (M') {\rotatebox{30}{\text{Picture language}}};
\node at (A1) {\rotatebox{60}{\text{Algebraic geometry}}};
\node at (A2) {\rotatebox{30}{\text{CFT 
}}};
\node at (A3) {\rotatebox{0}{\text{Fourier analysis}}};
\node at (A4) {\rotatebox{-30}{\text{Functional analysis}}};
\node at (A5) {\rotatebox{-60}{\text{Lattice gauge theory}}};
\node at (A6) {\rotatebox{-90}{\text{Knot theory}}};
\node at (A7) {\rotatebox{-120}{\text{QFT }}};
\node at (A8) {\rotatebox{-150}{\text{Quantum information}}};
\node at (A9) {\rotatebox{-180}{\text{Reflection positivity}}};
\node at (A10) {\rotatebox{-210}{\text{Representation theory}}};
\node at (A11) {\rotatebox{-240}{\text{Statistical physics}}};
\node at (A12) {\rotatebox{-270}{\text{Subfactor theory}}}; 
\end{scope}
\end{tikzpicture}
}
\end{equation*}

\section{Questions\label{Sect:Questions}}
We learned from our picture journey, that focusing on the picture language itself is very fruitful.  So here we collect a few questions for the future. 
The quon language provides an example, but we leave open the possibility of having many useful languages.

\subsection{Some Big Picture Questions for Birds}
\begin{enumerate}\footnotesize
\item{}  How far can one understand mathematical duality in terms of pictorial duality?   For example, to what extent can one understand further properties of Fourier duality or  mirror symmetry?

\item{}   Many modern picture  languages concern discrete combinatorial or topological data. A big question is: how to construct a continuum theory from those pictures?  Then one can ask  how one can understand continuous  symmetries, such as rotation invariance, in terms of picture language.  

\item{}  Many people have studied pictures from the point of view of topology and algebra.  
How far can one go  
to understand a different aspect: the analysis of pictures?

\item{}  Can one construct a CFT from a unitary MTC?

\end{enumerate}

\subsection{Some Technical Questions About  R}
\begin{enumerate}\footnotesize
\item{} Which family of mathematical concepts in {\bf R} is pictorial? 

\item{}  Given symmetries and related identities in {\bf R} which one would like to understand,  can one find pictures in {\bf L} that reflect these symmetries?

\item{}  Can one identify identities in {\bf R}  in terms of elementary operations on pictures in {\bf L}, such as by topological isotopy?

\item{}  Can we do computations for {\bf R}  in {\bf L} without using {\bf R}?  In this case we say that the language {\bf L} is complete.
Complete picture relations may require virtual concepts that do not make sense in {\bf R}.

\item{}  Ideally the real concepts in {\bf R} can be represented by simple pictures in {\bf L}.  But in case there are virtual pictures in {\bf L}, can one find another {\bf R} and {\bf S} such that the virtual pictures become real?  This situation leads to relations between different subjects.
\end{enumerate}

\subsection{Some Technical Questions About  L}
\begin{enumerate}\footnotesize
\item{}  How does one define a new picture language?  The main constraint is consistency.  One needs to introduce axioms compatible with pictorial intuition.   Sometimes one requires additional axioms that are motivated by requirements learned from~{\bf R}.

\item{}  How does one construct examples of a picture language, possibly using one of the following methods?
(i)~Find  {\bf L}~from~{\bf R}.
(ii) ~Derive new examples from known ones.
(iii)~Construct examples in an abstract way.

\item{}  How does one study pictures using axioms, as in Euclidean geometry?
One can thus ask questions about the language, based on its own interest.  

\item{}  After defining a pictorial language {\bf L} in an abstract way, how can one use the language to simulate interesting mathematics? For example, is there a CFT associated with the Haagerup subfactor~\cite{Haagerup}?
\end{enumerate}

\subsection{Some Other Questions}
\begin{enumerate}\footnotesize
\item{}
We have seen that the CNOT  gate in quantum information (illustrated here in the Quon language) has a pictorial interpretation.
\begin{align*}
\raisebox{-1cm}{
\scalebox{.8}{
\begin{tikzpicture}
\begin{scope}[xscale=.7,yscale=.7]
\begin{scope}[shift={(4,2)}]
\begin{scope}[shift={(0,0)},xscale=.5,yscale=.5]
\pgftransformcm{1}{0}{0}{1}{\pgfpoint{0}{0}}
\draw[blue!100,dashed] (0,0,-6)--++(0,0,9);
\draw[blue,dashed] (0,0,3)--++(0,0,6);
\draw[blue!100] (0,.8,-6)--++(0,0,9)--+(0,1,0) ;
\draw[blue!100] (1,0,-6)--++(0,0,9);
\draw[blue] (1,0,3)--++(0,0,6);
\draw[blue!100] (1,.8,-6)--++(0,0,10)--+(0,1,0);
\draw[blue] (0,.8,9)--++(0,0,-4)--+(0,1,0);
\draw[blue] (1,.8,9)--++(0,0,-3)--+(0,1,0);
\draw[white,WL] (0,1.8,3)--++(0,1,0)--+(-4,0,-4);
\draw[purple,dashed] (0,1.8,3)--++(0,1,0)--+(-4,0,-4);
\draw[white,WL] (0,1.8,5)--++(0,1,0)--+(-4,0,-4);
\draw[purple] (0,1.8,5)--++(0,1,0)--+(-4,0,-4);
\draw[white,WL] (1,1.8,4)--++(0,1,0)--+(6,0,6);
\draw[purple,dashed] (1,1.8,4)--++(0,1,0)--+(6,0,6);
\draw[white,WL] (1,1.8,6)--++(0,1,0)--+(6,0,6);
\draw[purple] (1,1.8,6)--++(0,1,0)--+(6,0,6);
\draw [white,WL] (-4,3.6,-1)--++(11,0,11);
\draw [purple] (-4,3.6,-1)--++(11,0,11);
\draw [white,WL] (-4,3.6,1)--++(11,0,11);
\draw [purple] (-4,3.6,1)--++(11,0,11);
\end{scope}
\end{scope}
\end{scope}
\end{tikzpicture}}}
\quad \quad \text{Quon\ CNOT}
\end{align*}
Another important gate in quantum information is the Toffoli gate.  Is the Toffoli gate topological?

\item{} In the Quon language, one illustrates a qudit by a 3D picture where the Z and X coordinate directions play a special role.
Is there a representation of a qudit by (3+1)D picture in which the 3 Frobenius algebras associated with X, Y, Z are represented by pictures in three orthogonal directions, and such that associativity  becomes 3D topological isotopy?

\item{}  Biamonte has posed other quantum information questions in~\cite{Biamonte}, including a pictorial understanding of the Gottesman-Knill theorem.

\item{}  Is there a picture language with  pictorial representations of differentiation, ordinary differential equations, and  partial differential equations?

\item{}  Feynman diagrams give pictures for Hermite polynomials for a given Gaussian.  Is there a similar understanding of all Gaussians, Fourier transforms, and their Hermite polynomials?

\item {} Does a pictorial lattice model have a continuum limit as an interesting quantum field theory?
If reflection  positivity is involved, can one preserve positivity at each step in the approximation and the construction of a limit?

\item{}  Kramers-Wannier duality allows one to compute the critical temperature for the 2D Ising model.  Is there a duality that allows one to compute the critical temperature of a pictorial statistical model or its field theory limit?
Is the limit a CFT at the critical temperature?
\end{enumerate}

\section{Acknowledgement}
We are grateful to Lusa Zheglova for the use of her drawing. 
We thank the Operator Algebra and Subfactor program at the Isaac Newton Mathematical Institute for hospitality. 
This research was supported in part by grants TRT0080 and TRT0159 from the Templeton Religion Trust.

\end{document}